\documentclass[12pt]{article}
\usepackage{amsmath,amsfonts,amssymb}
\usepackage[margin=1.2in]{geometry}

\newtheorem{proposition}{Proposition}[section]
\newtheorem{definition}{Definition}[section]

\newtheorem{theorem}{Theorem}

\usepackage{graphicx}
\usepackage[margin=6pt,skip=0pt,font=small,labelfont=bf,labelsep=period]{caption}
\usepackage[dvips]{hyperref}

\newcommand{\lmax}{\Lambda_{\max}}
\newcommand{\R}{\mathbb R}
\usepackage{soul}

\begin{document}

\title{Dynamics of Periodically-kicked Oscillators}

\author{Kevin K. Lin\thanks{Department of Mathematics and Program in
    Applied Mathematics, University of Arizona, AZ 85721, USA.  E-mail:
    klin@math.arizona.edu.  KL was supported in part by NSF Grant
    DMS-0907927.}\ \ \ and\ \ Lai-Sang Young\thanks{Courant Institute of
    Mathematical Sciences, New York University, NY 10012, USA.  E-mail:
    lsy@cims.nyu.edu.  LSY was supported in part by NSF Grant
    DMS-0600974.}}

%\date{\input{kickedosc.date}}
\date{\today}

\maketitle

\vskip -.1in

\centerline{\it Dedicated to Steve Smale on the occasion of his
80th birthday}

\vskip .3in

\begin{abstract}
  We review some recent results surrounding a general mechanism for
  producing chaotic behavior in periodically-kicked oscillators.  The
  key geometric ideas are illustrated via a simple linear shear model.
  \end{abstract}

\section*{Introduction}

This paper reviews some recent results on a topic with a considerable
history: the periodic forcing of limit cycles.  Some 80 years ago, van
der Pol and van der Mark observed that irregularities developed when
certain electrical circuits exhibiting stable oscillations were
periodically forced \cite{van der pol}.  Their work stimulated a number
of analytical studies; see {\it e.g.} \cite{CL, Ln, Li, H}.
Another classical example of driven oscillators is the FitzHugh-Nagumo
neuron model~\cite{F};
the response of this and other models of biological rhythms  to external
perturbations have been extensively studied (see {\em e.g.}
\cite{Winfree}).  As a topic of mathematical study, the dynamics of
forced oscillations is well motivated:
Oscillatory behavior are ubiquitous in physical, biological, and
engineered systems, and external forcing, whether artificially applied
or as a way to model forces not intrinsic to the system, is also commonplace.

In this article, we are not concerned with modeling specific physical
phenomena.  Instead, we consider a generic dynamical system with a limit
cycle, and seek to understand its qualitative behavior when the system
is periodically disturbed. To limit the scope of the problem, we
restrict ourselves to periodic {\it kicks}, or forcings that are turned
on for only short durations, leaving the limit cycle ample time to
restore itself during the relaxation period.  We are interested in
large-time behavior, particularly in questions of stability and chaos.
 
As we will show, one of the properties of the limit cycle that plays a key
role in determining whether the kicked system is stable
or chaotic is {\it shear}, by which we refer to the differential
in speed (or angular velocity) for orbits near the limit cycle. 
A central theme of this article is that under suitable conditions, the impact of
a kick can be substantially magnified by the underlying shear 
in the unforced system, leading to
the formation of horseshoes and strange attractors. This does not always 
happen, however: some limit cycles are more vulnerable, and some
types of kicks are more effective than others.
These ideas are discussed in \cite{WY2,WY3} and \cite{LY}, the material which
forms the basis of the present review. 

Even though we seek to provide insight into dynamical mechanisms that
operate under general conditions, we have found the ideas to be most
transparent in a very simple linear shear model, to which we will devote
a nontrivial part of the paper.  Sect.~1 introduces this example and
familiarizes the reader with the various parameters (including the one
which measures shear); it also reports on results of a numerical study
on Lyapunov exponents.  Sects.~2 and 3 are organized around explaining
these simulation results. Along the way, we take the opportunity to
review a number of mathematical ideas which clearly go beyond this one
example. Some of the rigorous results reviewed, notably those on SRB
measures for a relevant class of strange attractors \cite{WY1, WY4,
  WY5}, are recent developments. With the main ingredients of the linear
shear model and the relevant mathematical background in hand, we return to a discussion of general limit cycles in the
final section.

\section{Increasing Shear as a Route to Chaos}

%The aim of this section is to

This section introduces the main example we use in this review, and
acquaints the reader with the various parameters in the model and how
they impact the dynamics.  Of particular interest to us is the effect of
increasing {\it shear}.  Numerically computed Lyapunov exponents as
functions of shear are presented in Sect.~1.3. They will serve as a
focal point for some discussions to follow.

\subsection{Periodic kicking of linear shear flow}

%% As we explained in the Introduction, it is best to focus initially on
%% a specific model, namely the periodic kicking of a linear shear flow
%% with a hyperbolic limit cycle. This 2D model was studied rigorously
%% in~\cite{WY2,WY3} and numerically in~\cite{lin-young,zaslav}.  Though
%% exceedingly simple in appearance, it already exhibits rich and
%% complex dynamical behaviors.

Our main example is the periodic
kicking of a linear shear flow with a hyperbolic limit cycle. This 2D
model was studied rigorously in~\cite{WY2,WY3} and numerically
in~\cite{LY,zaslav}.  Though exceedingly simple in appearance, it
already exhibits rich and complex dynamical behaviors.

%% As noted in the Introduction, this is the simplest of a class of
%% systems that exhibit a rich and complex set of dynamical behaviors
%% associated with shear.

The model
%in question
is given by
\begin{equation}
\begin{array}{ccl}
\dot{\theta} &=& 1 + \sigma y\ ,\\[1ex]
\dot{y} &=& -\lambda y~+~A\cdot \sin(2\pi\theta)
\cdot\sum_{n=-\infty}^\infty \delta(t-n\tau)\ ,\\
\end{array}
\label{model1}
\end{equation}
where $(\theta, y) \in {\mathbb S}^1 \times {\mathbb R}$, 
${\mathbb S}^1 \equiv {\mathbb R}/{\mathbb Z}$, and $\sigma, \lambda, A$ 
and $\tau$ are constants with $\sigma, \lambda, \tau >0$. We will refer to 
Eq.~(\ref{model1}) with $A=0$ as the {\it unforced equation}, and the
term involving  $A$ as the {\it forcing} or the {\it kick}. Here $\delta$ is the usual
$\delta$-function, that is to say, the kicks occur instantaneously at times
$0,\tau, 2\tau, 3\tau, \dots$. 

More precisely, let $\Phi_t$ denote the flow corresponding
to the unforced equation. It is easy to see that for all $z \in {\mathbb S}^1 \times 
{\mathbb R}$, $\Phi_t(z)$ tends to the limit cycle $\gamma =\{y=0\}$ as $t \to \infty$.
The precise meaning of Eq.~(\ref{model1}) is as follows: Let 
$\kappa(\theta, y) = (\theta, y+A\sin(2\pi\theta))$ be the
kick map; it represents the action of the forcing term. 
Then assuming we start with a kick at time $0$, 
the time-$\tau$ map of the flow
generated by Eq.~(\ref{model1}) is $$\Psi_\tau = \Phi_\tau \circ \kappa,$$ and the evolution of the system is defined by iterating $\Psi_\tau$. 
We generally assume that
$\tau$ is not too small, so that during the relaxation period between kicks,
the flow $\Phi_t$
of the unforced equation ``restores" the system to some degree. 

The parameters of interest are:
\begin{align*}
\sigma &= \mbox{amount of shear,}\\
\lambda &= \mbox{rate of contraction to $\gamma$,}\\
A &= \mbox{amplitude of kicks, and}\\
\tau &= \mbox{time interval between kicks.}
\end{align*}
Our aim in the remainder of this section is to understand -- via
geometric reasoning and numerical simulations --
%% how the dynamical behavior of the system Eq.~(\ref{model1}) depends on
%% these parameters.  We seek to understand
the meanings of these quantities, and the roles they play in questions
of stability and chaos.  Our line of reasoning follows \cite{WY2} and
\cite{LY}.

%%%%%%%%%%%%%%%%%%%%%%%%%%
\subsection{Geometry of $\Psi_\tau$}

A simple way to gain intuition on the geometry of $\Psi_\tau$ is 
to study the $\Psi_\tau$-image of $\gamma$, the limit cycle of the unforced system.
We will do so by  freezing some of the parameters and varying others.

\bigskip
\noindent {\bf Effects of varying $\sigma, \lambda$, and $A$}

\smallskip
To begin with, let us freeze $\lambda,$ $A,$ and $\tau$. To fix ideas,
let us take $\lambda$ to be relatively small, so that the rate of
contraction is weak, and choose $\tau$ large enough that $e^{-\lambda
  \tau}$ is a nontrivial contraction.  This is when the effects of shear
are seen most clearly.  In Fig.~1, $\lambda = A=0.1$, and $\tau=10$.
Here $\Psi_\tau(0,0)=(0,0)$ because the limit cycle $\gamma$ has period
$1$ and $\tau$ is an integer multiple of this period; for non-integer $\tau$ the
picture is shifted horizontally.

Fig.~1(b) shows four images of $\gamma$ under $\Psi_\tau$ for increasing 
shear. The larger $\sigma$, the greater the 
%differential in instantaneous speeds of two points whose $y$-coordinates differ. 
difference in velocity
 between two points with different $y$ coordinates.
This applies in particular to the highest and lowest points
 in $\kappa(\gamma)$ in Fig.~1(a).
For $\sigma$ small enough, order in the $\theta$-direction is preserved,
{\it i.e.}, for $z_1 = (\theta_1,0)$ and $z_2=(\theta_1 + \varepsilon, 0)$,
$\Psi_\tau(z_1)$ will continue to have a slightly smaller $\theta$-coordinate than
$\Psi_\tau(z_2)$. As $\sigma$ increases, some points in $\gamma$ may 
``overtake" others, spoiling this order. As $\sigma$ gets larger still,
the total distances traveled in $\tau$ units of time vary even more, 
and a fold develops. 
This fold can be made arbitrarily large: we can make it wrap around the 
cylinder  as many times as we wish by taking $\sigma$ large enough.

%%%%%%%%%%%%%%%%%%%%%%%%%%%%%%
%% All 5 images are (i) scaled to 90% of the original, and (ii) squashed
%% VERTICALLY by another 30%.
\begin{figure}
  \begin{center}
    \begin{tabular}{c}
      \begin{tabular}{r}
        $y$ \\[1.4in]
      \end{tabular}
      \includegraphics*[scale=0.9,bb=14pt 60pt 3.1in 190pt,width=209.2pt,height=91pt]{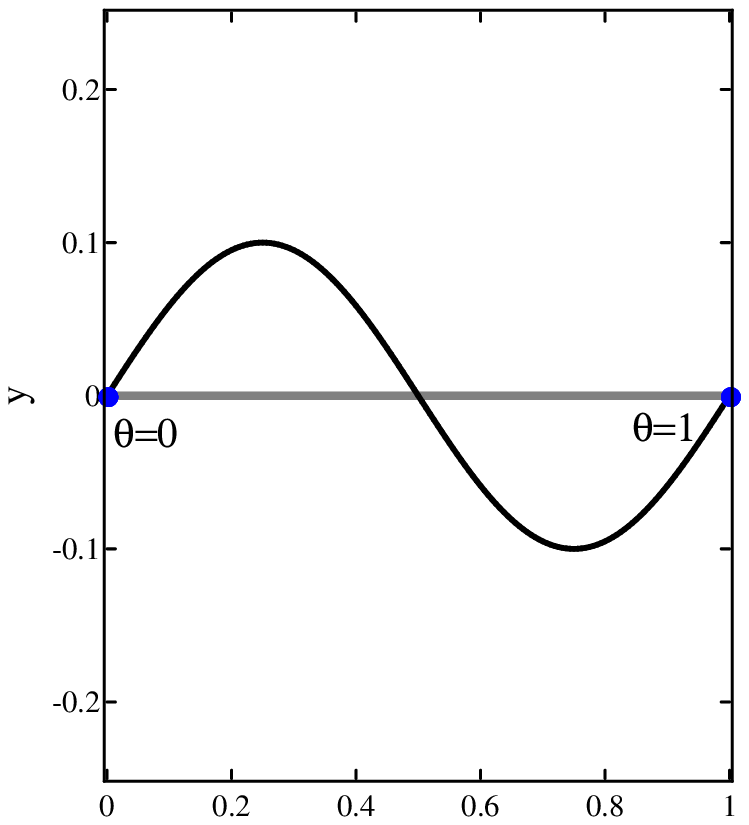} \\[-.7in]
      %\hspace{20pt}$\theta=0$\hspace{80pt}$\frac12$\hspace{80pt}1 \\[2ex]
    \end{tabular}

    {\small (a) The limit cycle $\gamma$ and its image $\kappa(\gamma)$
    after one kick}
    \\[0.25in]

    \begin{tabular}{lp{4ex}l}
      $\sigma=0.05$ && $\sigma=0.25$ \\[1ex]
      \includegraphics[scale=0.9,bb=89 -5 272 35,width=183pt,height=28pt]{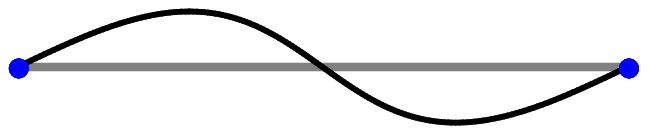} 
      &&\includegraphics[scale=0.9,bb=89 -5 272 35,width=183pt,height=28pt]{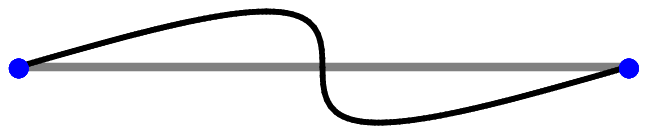} \\[4ex]
      $\sigma=0.5$ && $\sigma=1$ \\
      \includegraphics[scale=0.9,bb=89 -5 272 35,width=183pt,height=28pt]{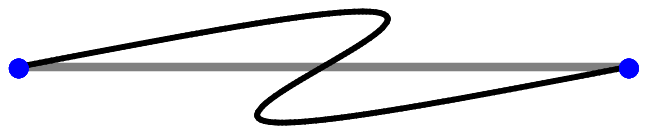} 
      &&\includegraphics[scale=0.9,bb=89 -5 272 35,width=183pt,height=28pt]{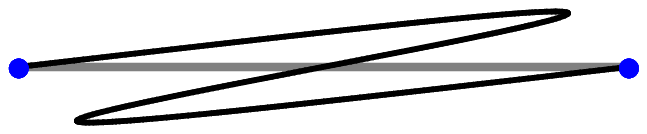} \\[4ex]
    \end{tabular}

    {\small (b) $\Psi_\tau(\gamma)$ for $\tau=10$}
  \end{center}

  \caption{Effect of increasing shear.  Here, $\lambda=A=0.1$.}
\end{figure}
%%%%%%%%%%%%%%%%%%%%%%%%%%%%%%%%

If we had fixed $\sigma$ instead, and increased $\lambda$ starting from
$\lambda=0.1$, the resulting sequence of pictures would be qualitatively
similar to Fig.~1(b) but in reverse order: The smallest $\lambda$ would
correspond to the bottom-right image in Fig.~1(b), and the largest
$\lambda$ to the top-left --- provided $\tau$ is scaled so that $\lambda
\tau$ remains constant.  This is because for $\lambda$ small,
$\kappa(\gamma)$ returns to $\gamma$ very slowly, giving the shear a
great deal of time to act, while for larger $\lambda$, $\kappa(\gamma)$
is brought back to $\gamma$ more quickly. Thus all else being equal, 
$\sigma$ and $\lambda$, {\it i.e.}, shear and damping, have opposite
  effects.
  
The consequence of varying $A$ while keeping the other parameters fixed is
easy to see: the stronger the kick, the greater the difference in
$y$-coordinate between the highest and lowest points in
$\kappa(\gamma)$, and the farther apart their $\theta$-coordinates will
be when flowed forward by $\Phi_\tau$.

What we learn from the sequence in snapshots in Fig.~1 is that
{\it $A$ acts in concert 
with $\sigma$ to promote fold creation, while $\lambda$ works against it.}

\bigskip
\noindent{\bf Formulas for $\Psi_\tau$}

\smallskip
Since the unforced equation is easy to solve, one can in fact write down
explicitly the formulas for $\Psi_\tau$. Let $(\theta_\tau, y_\tau)=\Psi_\tau(\theta_0, y_0)$. A simple
  computation gives
\begin{equation}
\label{time-t map}
\begin{array}{rcl}
\theta_\tau & = & \theta_0 + \tau + \frac{\sigma}{\lambda} \cdot 
[ y_0 + A \sin(2\pi\theta_0)] \cdot (1 - e^{-\lambda
  \tau})\qquad\mbox{(mod 1)}\ ,\\[12pt]
y_\tau & = & e^{-\lambda \tau} [y_0 + A \sin(2\pi\theta_0)]\ .
\end{array}
\end{equation}
The reader can easily check that Eq.~(\ref{time-t map}) is in agreement with
the intuition from earlier. 

Note the appearance of the ratio
$\frac{\sigma}{\lambda} A$, or rather 
$\frac{\sigma A}{\lambda}(1-e^{-\lambda \tau})$, in the nonlinear term in the equation for
$\theta_\tau$: the size of this term is a measure of the tendency for a fold to
develop in $\Psi_\tau(\gamma)$. As is well known to be the case, stretch-and-fold
is a standard mechanism for producing chaos. One can, therefore, think of 
the ratio 
$$
\frac{\sigma}{\lambda} A \ = \ \frac{\rm shear}{\rm contraction} \cdot {\rm kick \
amplitude}
$$ 
as the key to determining whether 
the system is chaotic, provided $\lambda \tau$ is large enough
that the factor $1-e^{-\lambda \tau}$ is not far from $1$. 

\bigskip
\noindent {\bf Trapping region and attractor }

\smallskip

 From the above, it is evident that much of the action takes
place in a neighborhood around $\gamma.$ Let ${U}=\{ |y| \le
A(e^{\lambda \tau} -1)^{-1}\}$, so that $\Psi_\tau({U}) \subset {U}$,
and define
$$\Gamma = \cap_{n \ge 0}\Psi_\tau^n({U})$$ 
to be the {\it attractor} for the system Eq.~(\ref{model1}). 
The basin of attraction of $\Gamma$ is the entire cylinder
$\mathbb S^1 \times\mathbb R$, since every orbit will eventually enter
$U$. This usage of the word ``attractor"  
implies no knowledge of dynamical indecomposability
(a condition required by some authors).

%%%%%%%%%%%%%%%%%%%%%%%%%%%%%%%
\subsection{Lyapunov exponents}

Another measure of chaos is orbital instability, or the speed at which
nearby orbits diverge. In this subsection, we focus on the larger of the
two Lyapunov exponents of $\Psi_\tau$, defined to be
$$
\Lambda_{\max}(z) \ = \ \lim_{n\to \infty} \frac{1}{n} \log \|D\Psi^n_\tau(z)\| \ = \ 
\sup_{v \ne 0} \ \lim_{n\to \infty} \frac{1}{n} \log \|D\Psi^n_\tau(z) \cdot v\|\ .
$$
Leaving technical considerations for later (see Sect.~3.1),  
we compute numerically $\Lambda_{\max}$ for the systems 
in question, sampling at various points $z \in {U}$. 
Notice that $\Lambda_{\max}$ measures the rate of divergence of nearby
orbits {\it per kick}, not per unit time. 

Each of the plots in Fig.~2 shows $\Lambda_{\max}$ as a function of
$\tau$ for the values of $\sigma, \lambda$ and $A$ specified. In all six
plots, we have fixed $\lambda=A=0.1$, while $\sigma$ varies from plot to
plot.  The first 4 values of $\sigma$ used in Fig.~2 are the same as
those used in Fig.~1(b).  In each plot, 10 randomly chosen initial
conditions are used, the largest and smallest computed values of
$\Lambda_{\max}$ are discarded, and the largest and smallest of the
remaining 8 values are shown, the smallest as a solid black dot, and the
largest, if visibly different than the smallest, as an open square. The
plots show $\tau \in [5,15]$. We give some idea of the rates of
contraction and sizes of the trapping regions ${U}$ for these
parameters: at $\tau=5$, $e^{-\lambda \tau} = e^{-0.5} \approx 0.61$,
and ${U} = \{|y| \lesssim 0.15\}$; at $\tau=15$, $e^{-\lambda \tau} =
e^{-1.5} \approx 0.22$, and ${U}=\{|y|\lesssim 0.03\}$.

%%%%%%%%%%%%%%%%%%%%%%%%%%%%%%%%
\begin{figure}
  \begin{center}
    \begin{tabular}{lp{0.25in}l}
      (a) $\sigma=0.05$ && (b) $\sigma=0.25$ \\[6pt]
      \includegraphics[bb=0in 0in 3.75in 3in,scale=0.6666]{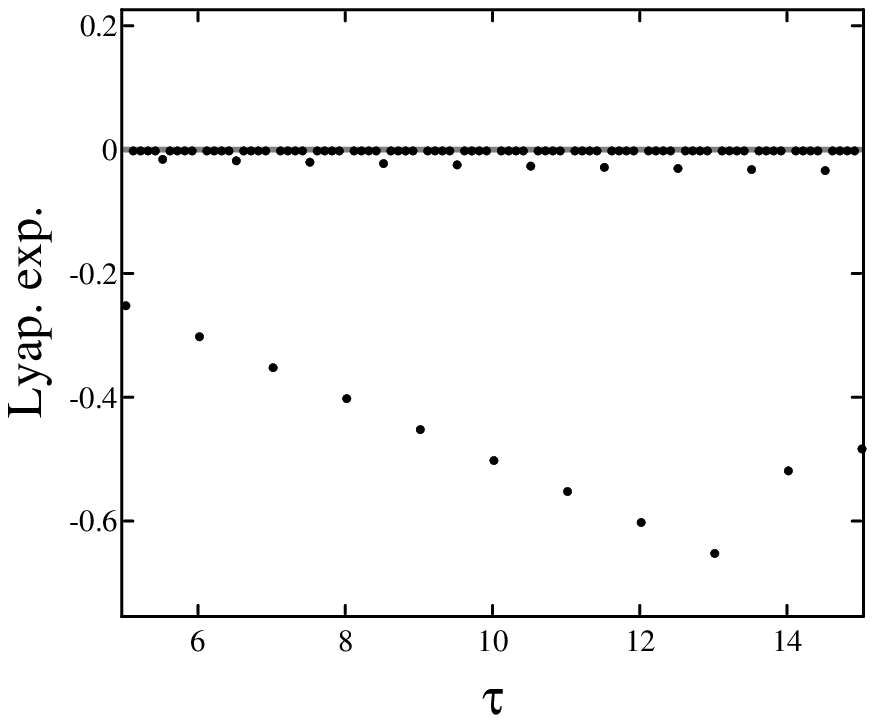} &&
      \includegraphics*[bb=0.5in 0in 3.75in 3in,scale=0.6666]{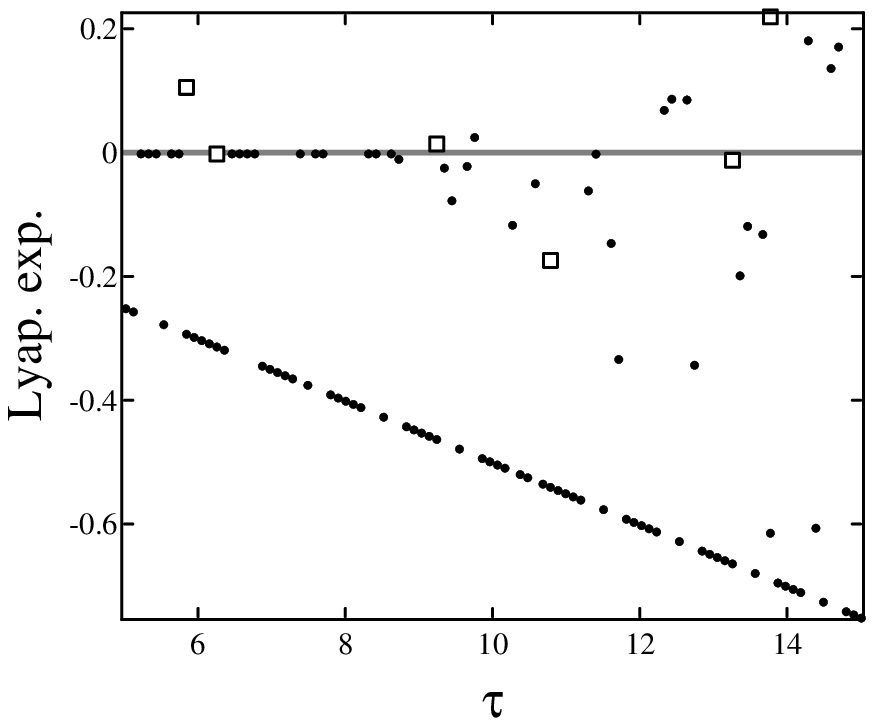} \\[12pt]
      (c) $\sigma=0.5$ && (d) $\sigma=1$ \\[6pt]
      \includegraphics[bb=0in 0in 3.75in 3in,scale=0.6666]{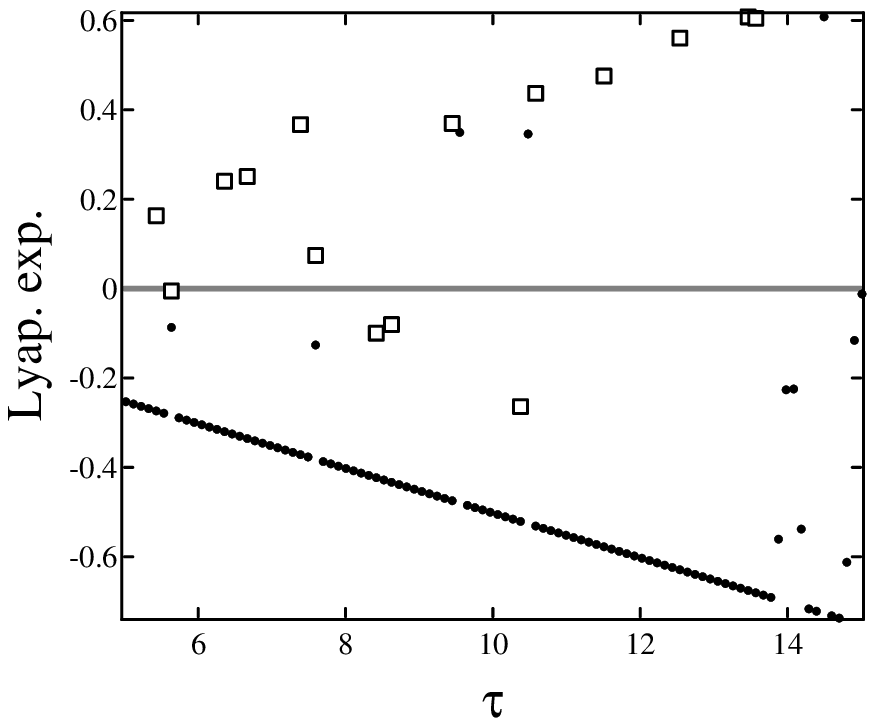} &&
      \includegraphics*[bb=0.5in 0in 3.75in 3in,scale=0.6666]{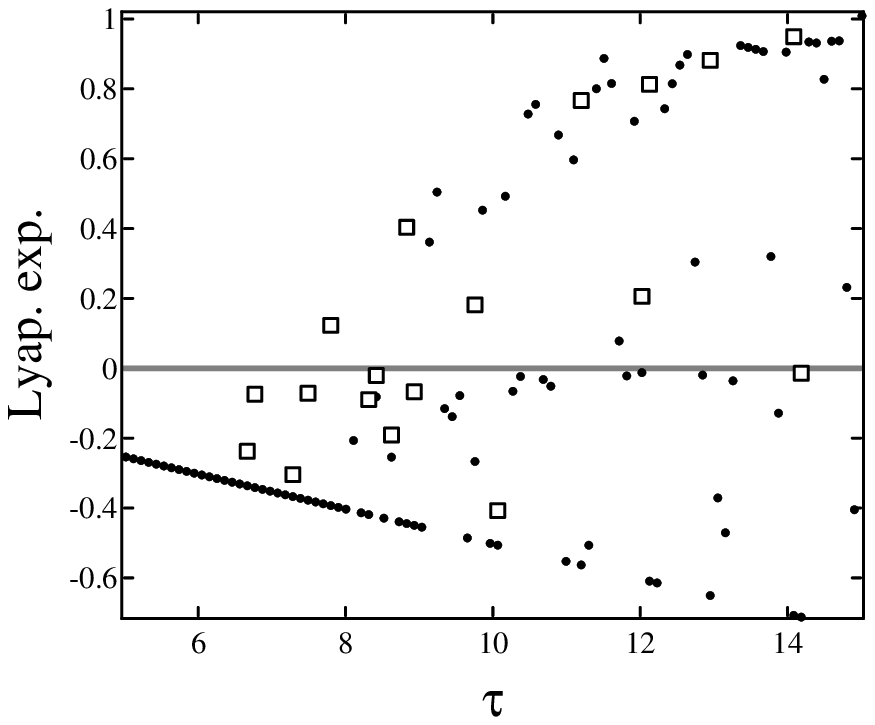} \\[12pt]
      (e) $\sigma=2$ && (f) $\sigma=4$ \\[6pt]
      \includegraphics[bb=0in 0in 3.75in 3in,scale=0.6666]{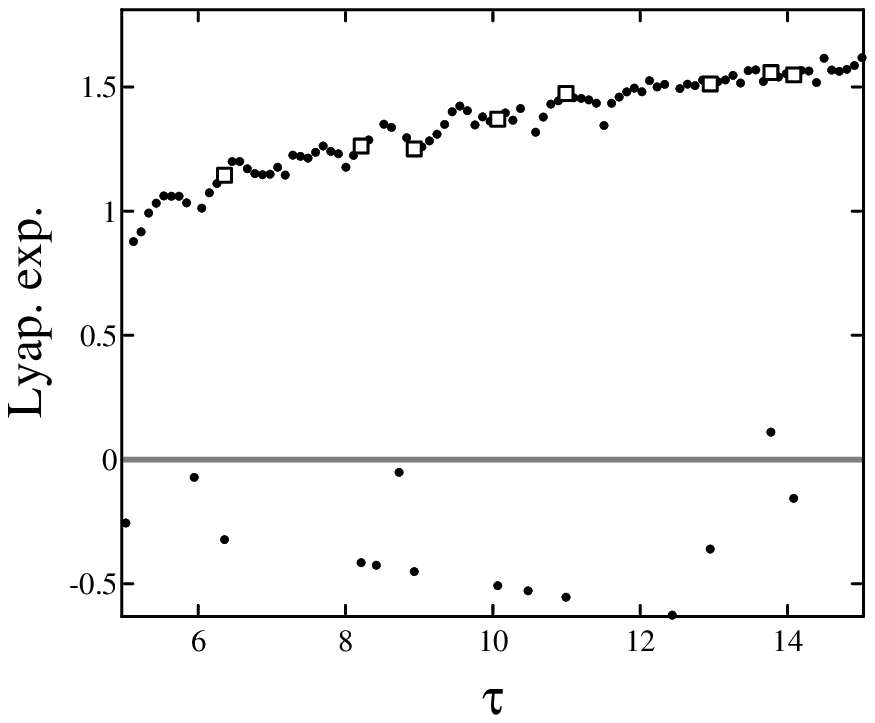} &&
      \includegraphics*[bb=0.5in 0in 3.75in 3in,scale=0.6666]{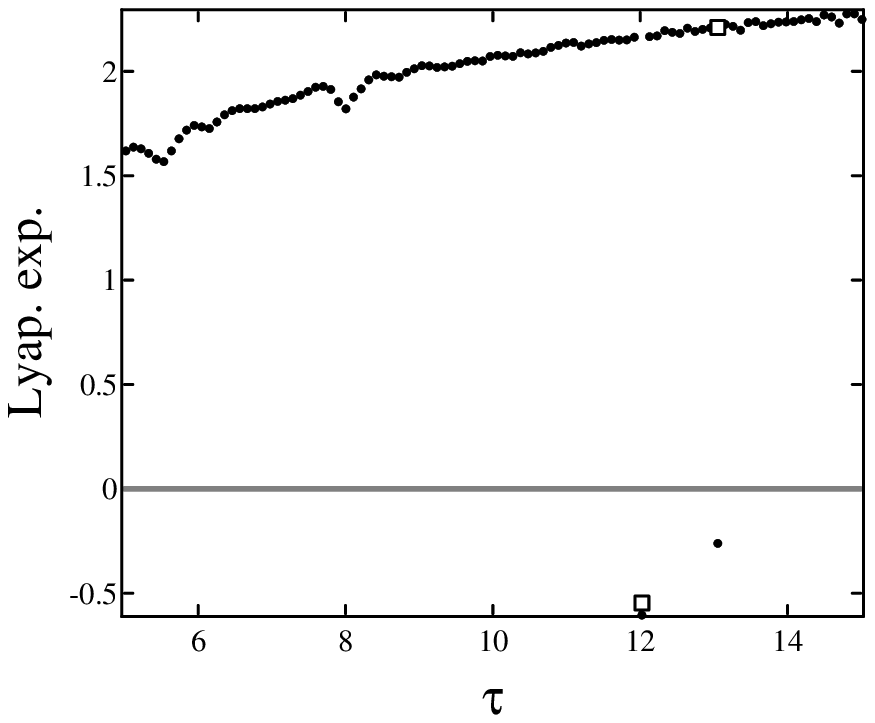} \\
    \end{tabular}
  \end{center}
  \caption{Lyapunov exponents $\Lambda_{\max}$ as functions of kick
    period $\tau$ for increasing shear $\sigma.$ The other parameters
    are $\lambda = A = 0.1.$ For each $(\sigma,\tau)$, we simulate 10
    orbits with random initial conditions, iterating each for $4\times
    10^6$ steps.  We then drop the outliers and plot the remaining
    estimates, as described in the text.}

    %% We drop the highest and lowest estimates and plot the remainder.
    %% For those $(\sigma,\tau)$ where not all orbits give essentially
    %% the same estimate, we have marked the upper estimate by a {\tiny
    %% $\square.$}
\end{figure}
%%%%%%%%%%%%%%%%%%%%%%%%%%%%%%%%

\bigskip
\noindent {\bf Observations from simulation results}

\smallskip
The discussion in Sect.~1.2 suggests that as shear is increased
with other parameters fixed, the system is likely to get increasingly chaotic.
This may lead us to expect $\Lambda_{\max}$ to increase monotonically
with $\sigma$. As one can see, that may be 
correct as an overall trend, but the situation is somewhat more complicated: 

In the two low-shear regimes, namely $\sigma = 0.05$ and $0.25$, 
with a few exceptions the computed values of 
$\Lambda_{\max}$ are either zero or negative, with a majority of them at or
very near zero for $\sigma = 0.05$ and becoming considerably more
negative at $\sigma = 0.25$. That is to say, 
$\Lambda_{\max}$ {\it decreases} as $\sigma$ increases.
Notice that $\Lambda_{\max}<0$ means 
the trajectory tends to a {\it sink}, {\it i.e.}, a stable fixed point or periodic orbit. 

Increasing shear, we see in the middle row of Fig.~2 that at first
sinks dominate the landscape, giving way to more instances of
positive Lyapunov exponents, {\it i.e.},
chaotic behavior, as $\sigma$ increases. At $\sigma = 1$, 
the picture is very mixed, with $\Lambda_{\max}$
fluctuating wildly between positive and negative values as $\tau$ varies.
Notice also the nontrivial number of open squares, telling us that 
these parameters often support more than one type of dynamical behavior.

In the two higher-shear regimes, $\sigma = 2$ and $4$, 
$\Lambda_{\max}$ becomes more solidly positive, though 
occasional sinks are still observed. The route has been a messy one, 
but one could say that the transition to chaos is complete. 

As to the dependence on $\tau$, it appears that other things being equal,
longer
relaxation times between kicks allow the dynamical phenomenon
in effect to play out more completely: regardless of the sign of
$\Lambda_{\max}$, its magnitude increases with $\tau$ in each of 
the plots.

Finally, it is important to remember that the limit cycles used to produce
the results in Fig.~2 are {\it weakly attracting}, making them more
vulnerable to the effects of shear. Strongly attracting limit cycles are more robust, 
and larger kicks and/or shear will be needed to produce chaos.

\bigskip
{\bf In Sects. 2 and 3 we review some rigorous theory that supports
  the numerically computed values of $\Lambda_{\max}$ shown.} To avoid
technical assumptions, we will focus on the model in Sect.~1.1, leaving
generalizations to Sect. 4. As the reader will see, the mathematical
ideas go considerably beyond this one example.  On the other hand, {\it
  even for this simple model, state-of-the-art understanding is
  incomplete}. In the next two sections, we will vary $\sigma, \lambda,
A$ and $\tau$, and show that there are regions in the
parameter space for which a clear description of the dynamics is
available, and larger regions on which there is partial understanding.

%%%%%%%%%%%%%%%%%%%%%%%%%%%%%%%%%%%%%%%
%%%%%%%%%%%%%%%%%%%%%%%%%%%%%%%%%%%%%%%%
\section{Geometric Structures}

To analyze a dynamical system, it is often useful to begin by
identifying its most prominent structures, those that are a significant
part of the landscape.  Even when they do not tell the whole story,
these structures will serve as points of reference from which to explore
the phase space.  This section describes structures of this type for the
systems defined by Eq.~(\ref{model1}).

\subsection{Persistence of limit cycles at very low shear}

\begin{proposition}[\cite{WY2}] Given $\lambda \tau>0$, the following
hold for $\frac{\sigma}{\lambda}$ and $A$ sufficiently small:

(a) the attractor $\Gamma $  is a smooth, closed invariant curve near $\gamma$;

(b) every $\zeta \in {\mathbb S}^1 \times {\mathbb R}$ lies in the strong
stable curve $W_{\Psi_\tau}^{ss}(z)$ for some $z \in \Gamma$. 
\label{prop:2.1}
\end{proposition}

Here $W_{\Psi_\tau}^{ss}(z)=\{\zeta \in {\mathbb S}^1 \times {\mathbb R}
: \limsup_{n \to \infty} \frac{1}{n} \log d(\Psi_\tau^n(z), \Psi_\tau^n(\zeta)) \le 
-\lambda' \tau\}$ where $\lambda'$ is a constant $> \frac12 \lambda$.
Proposition~\ref{prop:2.1} follows from standard arguments in
stable and center manifolds theory; see {\em e.g.} \cite{HPS}.
 The idea is simple: From Eq.~(\ref{time-t map}), one
obtains
\begin{equation} 
\label{derivative}
D\Psi_\tau(\theta, y) =
\begin{pmatrix}
1+2\pi \frac{\sigma}{\lambda}A \cos(2\pi \theta)(1-e^{-\lambda \tau})~~~ & 
\frac{\sigma}{\lambda}(1-e^{-\lambda \tau}) \\[12pt]
e^{-\lambda \tau} 2\pi A \cos(2\pi \theta) & e^{-\lambda \tau}
\end{pmatrix}\ .
\end{equation}
Since invariant cones depending on $\frac{\sigma}{\lambda}$ and $\lambda \tau$
clearly exist when $A=0$, they will persist when $\frac{\sigma}{\lambda}$ and
$A$ are small enough. 

When $\Gamma$ is a smooth invariant curve, the dynamics 
on $\Gamma$ is given by the theory of circle diffeomorphisms. The situation for a smooth one-parameter family of circle diffeomorphisms  $\{f_\omega\}$ can be summarized as follows (see {\it e.g.} \cite{GH}):
Let $\rho(f_\omega)$ denote the rotation number of $f_\omega$. 
Then $\omega \mapsto \rho(f_\omega)$ is a {\it devil's staircase},
the flat parts corresponding to intervals of $\omega$ on which the rotation
number is rational. Moreover, the set of $\omega$ for which $\rho(f_\omega)$ is rational is typically 
open and dense, while the set of $\omega$ for which $\rho(f_\omega)$ 
is irrational has positive Lebesgue measure. 
When $\rho(f_\omega) \in {\mathbb Q}$, $f_\omega$ typically
has a finite number of periodic sinks and sources alternating on the circle;
these aside, every orbit converges to a periodic sink. When $\rho(f_\omega) \not \in {\mathbb Q}$, $f_\omega$ is topologically conjugate to an irrational rotation. 

The ideas above capture the
spirit of the dynamics when shear is small enough: 
Suppose $\Lambda_{\max}$ is computed
using an initial condition $\zeta \in U$, and $\zeta \in W^{ss}(z)$ for 
$z \in \Gamma$. Then $\Lambda_{\max}(\zeta) =
\Lambda_{\max}(z)$, and from the discussion above, the latter is either strictly
negative or zero depending on whether $\rho(\Psi_\tau|_\Gamma)$ is rational or
irrational. 

\bigskip
\noindent {\bf Breaking of invariant curves}

\smallskip
Now if we fix $\lambda$ and $ A$, and increase the shear $\sigma$
as is done in Fig.~2, the invariant cones -- and the 
invariant curve itself -- will break. 

Here is how it happens in this model for integer values of $\tau$:
For $\tau \in {\mathbb Z}^+$, $(\theta, y)=(\frac12, 0)$ is a fixed point of 
$\Psi_\tau$, and a simple computation shows that as $\sigma$ increases from $0$, 
the larger eigenvalue of this fixed point decreases from $1$ to 
$e^{-\frac12 \lambda \tau} = \sqrt{\det(D\Psi_\tau)}$, at which time 
the eigenvalues turn complex.
No invariant curve can exist after that. Geometrically, one can think of
the breaking of the invariant curve as being due to too much ``rotation" 
or ``twist" at this fixed point. 

Taking this observation a step further, one notes from Eq.~(\ref{derivative}) that 
the rotational action of $D\Psi_\tau(\theta, y)$ is strongest  at 
$\theta = \frac12$, where $\cos(2\pi\theta)=-1$. This suggests that for fixed
$\sigma, \lambda$ and $A$, invariant curves are the most vulnerable for integer values of $\tau$, where this strongest rotation occurs at a fixed point. 

\bigskip
\noindent {\bf Interpreting Figs.~2(a) and (b)}

\smallskip
At $\sigma = 0.05$, a majority of the $\Lambda_{\max}$
values computed are at or very slightly below zero. This is consistent with
the existence of an invariant curve for those $\tau$. One checks easily that 
for integer values of $\tau \le 13$, $\Lambda_{\max}=-\frac12 \lambda \tau$ 
and the eigenvalues are complex. These are the first places 
where the invariant curve is broken as predicted.

In the plot for $\sigma=0.25$, without pretending to account for
all data points, it looks as though many are dropping off the 
$\Lambda_{\max}=0$ line to join the $\Lambda_{\max}=-\frac12 \lambda \tau$
line. The only
holdouts for $\Lambda_{\max}=0$ occur for smaller $\tau$ where, as noted
earlier, shear has not had enough time to act. 

%%%%%%%%%%%%%%%%%%%%%%%%%%%%%%%%%%
\subsection{Increasing shear: horseshoes and sinks}

\noindent {\bf At first, mostly sinks}

\smallskip
Fig.~2(c),(d) suggest that at $\sigma = 0.5$, a sink with complex conjugate eigenvalues dominates the scene for much of the range of $\tau$ considered,
and the same is true at $\sigma=1$ for smaller values of $\tau$.

For $\tau \in {\mathbb Z}^+$, this again is easily checked. The ``twist" at 
$\theta = \frac12$ is also eminently visible in the last three pictures of
$\Psi_\tau(\gamma)$ in Fig.~1(b). With a little bit of work, one can settle
these questions rigorously, but an {\it a priori} fact that makes plausible
the extension of this sink to non-integer values of $\tau$ is that 
fixed point sinks with complex conjugate eigenvalues 
cannot disappear suddenly as parameters are varied: a bifurcation 
can occur only when these eigenvalues become real, {\it i.e.}, 
$\pm e^{-\frac12 \lambda \tau}$, (and a fixed point can vanish only when
one of its eigenvalues is equal to $1$). 

Finally, we remark that even though the sinks above clearly exert
nontrivial influence on the dynamics, other structures (competing
sinks, invariant sets etc.) may be present. The many open
squares in Fig.~2(c) suggest that for these parameters, trajectories in
different regions of the phase space have distinct futures.

%%%%%%%%%%%%%%%%%%%%%%%%%%%%%%%%
\begin{figure}
  \begin{center}
    \includegraphics[bb=0in 0in 5in 2in]{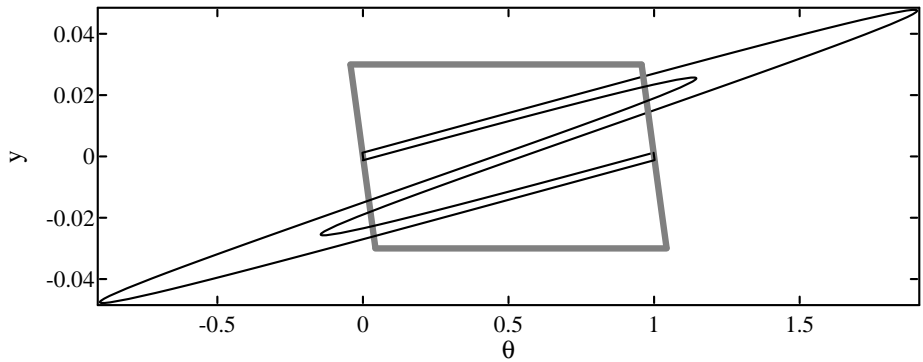} \\
    {\small (a) Formation of a horseshoe} \\[16pt]

    \begin{tabular}{c}
      $y=0$ \\[3in]
    \end{tabular}%
    %% BUG: If I generate a 5x3" eps file and print the resulting PDF on
    %% Macs using either Preview or Skim, the hairlines come out with
    %% nonzero width.  This appears to be a problem with Apple's PDF
    %% print driver -- the problem goes away with Acrobat.  If one
    %% generates a ps file directly, or re-convert the pdf to ps and
    %% then print THAT, the problem also goes away.
    %% Work-around: Make a big picture and scaling it down seems to fix
    %% it...
    \includegraphics[bb=0in 0in 10in 6in,scale=0.5]{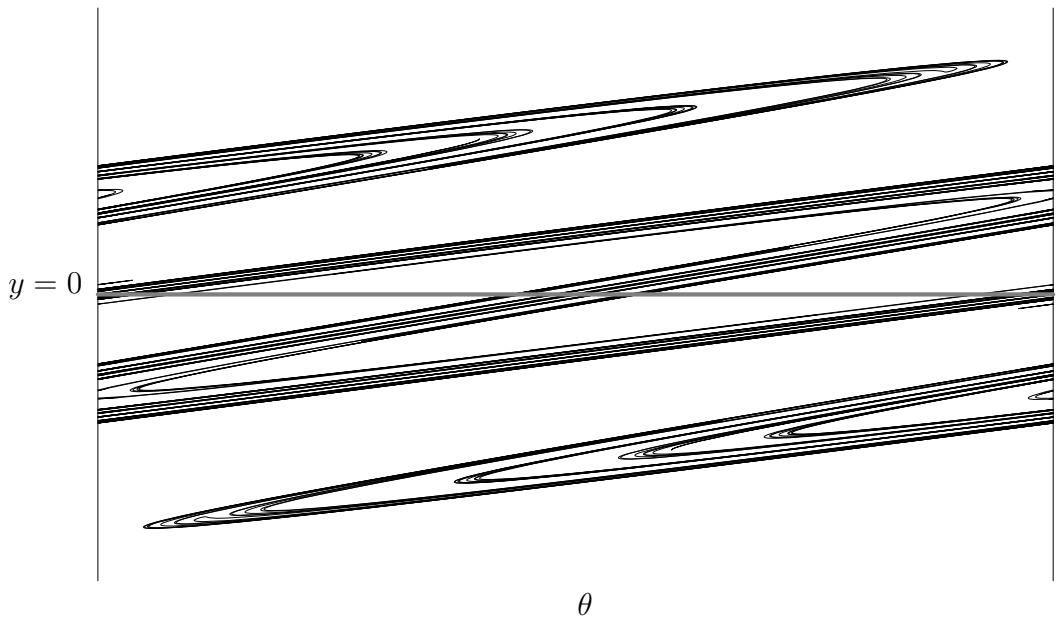} \\[-1.5in]
    ~~~~~~~~~~~~$\theta$ \\[3ex]
    
    {\small (b) The attractor $\Gamma$}
  \end{center}
  \caption{The attractor $\Gamma$ and a horseshoe in it.  Panel (a)
    illustrates the formation of a horseshoe: Shown are a box $R$ (thick
    gray lines) and its image $\Psi_\tau(R)$ (thinner black curves).
    The vertical boundaries of $R$ are the stable manifolds of the fixed
    point at $(0,0);$ these are mapped into themselves by $\Psi_\tau~.$
    Panel (b) shows a picture of the attractor $\Gamma.$
    %% The latter is computed by applying $\Psi^4_\tau$ to the line
    %% segment $(-\frac1{10},\frac1{10})\times\{0\}.$
    The parameters are $\lambda=A=0.1$, $\sigma=2,$ and $\tau=10$ for
    both plots. }
\end{figure}
%%%%%%%%%%%%%%%%%%%%%%%%%%%%%%

\bigskip
\noindent {\bf Smale's horseshoes}

\smallskip
Horseshoes are likely present starting from $\sigma$ somewhere between $0.5$
and $1$ for $\tau$ large enough. An example of an easily recognizable horseshoe for $\sigma = 2$ and $\tau=10$ is shown in Fig.~3(a).  
The larger $\sigma$, the easier it is to give examples. 

Fig.~3(a) illustrates how proofs of horseshoes or uniformly hyperbolic invariant
sets are often done: 
One first ``spots" a horseshoe with one's eyes, namely one or more boxes that map
across themselves in a characteristic way, and then proves that 
the set of points that remain in these boxes forever has the required splitting into expanding and contracting directions. When $\frac{\sigma}{\lambda}$ is large, 
$D\Psi_\tau$ expands strongly in the $\theta$-direction for most values of $\theta$;
contraction is guaranteed since $\det(D\Psi_\tau)=e^{-\lambda \tau}< 1$.

\begin{proposition}[\cite{WY2}] Given $\lambda \tau$, $\Psi_\tau$ has a horseshoe if
$\frac{\sigma}{\lambda}A$ is sufficiently large.
\end{proposition}

The presence of horseshoes is sometimes equated with dynamical
complexity or chaos in the literature, and that is entirely justified
insofar as one refers to the {\it existence} of chaotic orbits. One must
not confuse the existence of these orbits with chaotic behavior starting
from ``most" or ``typical" initial conditions, however: A system can
have a horseshoe (which attracts a Lebesgue measure zero set), and have
all other points in the phase space tending to a sink.  Or, the
horseshoe can be part of a ``strange attractor", with
$\Lambda_{\max}>0$.  The presence of a horseshoe alone does not tell us
which of these scenarios will prevail.  We will say more about
  strange attractors {\it versus} sinks in Sect.~3.1.
    Suffice it to observe here that horseshoes clearly exist
for most of the parameters in Fig.~2(d)-(f), and $\Lambda_{\max}$ is
sometimes positive and sometimes negative.

%% The presence of a horseshoe alone does not tell us which of these
%% scenarios will prevail.  This discussion is continued in Sect.~3.1.
%% Suffice it to observe here that horseshoes clearly exist for most of
%% the parameters in Fig.~2(d)-(f), and $\Lambda_{\max}$ is sometimes
%% positive and sometimes negative.

\bigskip
\noindent {\bf Sinks from homoclinic tangencies}

\smallskip
For larger shear, the attractor can be quite complicated; see Fig.~3(b). 
Yet in Figs.~2(d), (e), and even (f), sinks can also occur as noted.

The following is purely theoretical, in the sense that we do not know
exactly where the sinks are in these specific systems, but it is a
general fact in two dimensions that near homoclinic tangencies of
dissipative saddle fixed points (``dissipative" means $|\det(Df)|<1$),
sinks form easily, meaning one can perturb the map and find one near
such a tangency; see \cite{N}.  Furthermore, tangencies persist once the
stable and unstable manifolds of a horseshoe are shown to meet
tangentially somewhere. While no results have been proved for this
particular model, the ``turns" made by unstable manifolds (see
Fig.~3(b)) suggest the abundance of opportunities for such tangencies.

%%%%%%%%%%%%%%%%%%%%%%%%%%%%%%%%
%%%%%%%%%%%%%%%%%%%%%%%%%%%%%%%%
\section{A theory of strange attractors}

In this section we focus on the case of {\it positive} Lyapunov
exponents, having discussed negative and zero values of $\Lambda_{\max}$
in Section 2.  A combination of geometric and statistical ideas will be
used. Since these developments are more recent, we think it may be
useful to include more background information: Sect.~3.1 discusses SRB
measures for general chaotic systems. Sect.~3.2 surveys some recent work
on a class of strange attractors called {\it rank-one attractors}.  In
certain parameter ranges, the attractors $\Gamma$ in our kicked
oscillator systems are of this type.  In Sect.~3.3, we explain how the
general results reviewed in Sects.~3.1 and 3.2 are applied to
Eq.~(\ref{model1}).

\subsection{SRB measures}

The setting of this subsection is as follows:
Let $M$ be a Riemannian manifold or simply $\mathbb R^n$. We consider 
an open set ${U} \subset M$ 
with compact closure, and let $f$ be a $C^2$ embedding of ${U}$ into itself
with $\overline{f({U})} \subset {U}$. 
We will refer to $\Gamma = \cap_{n \ge 0} f^n({U})$ as the {\it attractor} and 
${U}$ as its {\it basin of attraction}. Though not a formal assumption, we have
in mind here situations that are ``chaotic"; in particular, $\Gamma$ is 
more complicated than an attracting periodic orbit. 

We will adopt the viewpoint that {\it observable events are represented
by positive Lebesgue measure sets,}  and are
 interested in invariant measures that reflect the 
properties of Lebesgue measure, which we denote by $m$.
For chaotic systems, the only invariant measures
known to have this property are SRB measures.
(The terms Lebesgue and Riemannian measures will be used interchangeably
in this article.) 

\begin{definition} An $f$-invariant Borel probability measure $\mu$
is called an {\rm SRB measure}~if

\smallskip
(a) $\Lambda_{\max}>0 \ \mu$-a.e., and

\smallskip
(b) the conditional measures of $\mu$ on local unstable manifolds
  have densities with 
  
\quad respect to the Riemannian measures on these manifolds.
\end{definition}

%\begin{definition} An $f$-invariant Borel probability measure $\mu$
%is called an {\rm SRB measure}~if
%\begin{enumerate}
%\item[(a)] $\Lambda_{\max}>0 \ \mu$-a.e., and
%\item[(b)] the conditional measures of $\mu$ on local unstable manifolds
%  have densities with respect to the Riemannian measures on these
%  manifolds.
%\end{enumerate}
%\end{definition}

Recall from the Multiplicative Ergodic Theorem~\cite{O} that Lyapunov
exponents, in particular $\Lambda_{\max}$, are defined $\mu$-a.e., so
(a) makes sense; in general, these quantities may vary from point to
point. The meaning of (b) can be understood as follows.  For an
invariant measure to reflect the properties of $m$, it is simplest if it
has a density with respect to $m$, but that is generally not possible
for attractors: All invariant measures in ${U}$ must live on $\Gamma$,
and if $f$ is volume decreasing, which is often the case near an
attractor, then $m(\Gamma)=0$.  The idea of SRB measures is that if
$\mu$ cannot have a density, then the next best thing is for it to have
a density in unstable directions, the intuition being that the
stretching of phase space in these directions leads to a smoothing of
distributions.

The main result on SRB measures is summarized in the next Proposition,
followed by a sketch of its proof.
The ideas in the proof will explain how SRB measures, which are themselves singular,
are related to Lebesgue measure.  Recall that for an ergodic measure
$\mu$, $\Lambda_{\max}$ is constant $\mu$-a.e.  We will denote this
number by $\Lambda_\mu$.

%% Its proof, a sketch of which will follow, explains how SRB measures,
%% which are themselves singular, are related to Lebesgue
%% measure. Recall that for an ergodic measure $\mu$, $\Lambda_{\max}$
%% is constant $\mu$-a.e.  We will denote this number by $\Lambda_\mu$.

\begin{proposition} Let $(f,\mu)$ be an ergodic SRB measure 
with no zero Lyapunov exponents.  Then there is a 
set ${V} \subset {U}$ with $m(V)>0$ such that the following hold for every $y \in
{V}$:

\smallskip
(i) $\Lambda_{\max}(y) = \Lambda_\mu$; and

\smallskip
(ii) $\frac{1}{n} \sum_{i=0}^{n-1} \varphi(f^iy)\to \int \varphi
  d\mu\ $ for every continuous observable $\varphi: {U} \to {\mathbb
    R}.$
\end{proposition}
%

%\begin{proposition} Let $(f,\mu)$ be an ergodic SRB measure 
%with no zero Lyapunov exponents.  Then there is a 
%set ${V} \subset {U}$ with $m(V)>0$ such that the following hold for every $y \in
%{V}$:
%\begin{enumerate}

%\item[(i)] $\Lambda_{\max}(y) = \Lambda_\mu$; and

%\item[(ii)] $\frac{1}{n} \sum_{i=0}^{n-1} \varphi(f^iy)\to \int \varphi
%  d\mu\ $ for every continuous observable $\varphi: {U} \to {\mathbb
%    R}.$

%\end{enumerate}
%\end{proposition}

The idea of the proof is as follows: Let $\gamma$ be a piece of local
unstable manifold, and let $m_\gamma$ be the Riemannian measure on
$\gamma$. By property (b) of Definition 3.1, we may assume
$m_\gamma$-a.e. $x \in \gamma$ is ``typical" with respect to $\mu$.  In
particular, it has properties (i) and (ii) in the Proposition.  Let
$W^s(x)$ be the stable manifold through $x$.  Properties (i) and (ii)
for $y\in W^s(x)$ follow from the corresponding properties for $x$ because
$d(f^nx, f^ny) \to 0$ exponentially as $n \to \infty$.  It remains to
show that the set of points $y$ that are connected to $\mu$-typical
points as above has positive $m$-measure, and that is true by the
absolute continuity of the stable foliation \cite{PS}. 

\medskip
A little bit of history:  SRB measures were invented by Sinai,
Ruelle and Bowen in the 1970s, when they constructed for every attractor 
satisfying Smale's Axiom A \cite{Sm} a special invariant measure with the properties 
in Definition 3.1 (\cite{Sn, R1, BR}).\footnote{Sinai treated 
first the case of Anosov systems; his results were shortly thereafter extended 
to Axiom A attractors (which are more general) by first Ruelle and then 
Ruelle and Bowen.} This special invariant measure has a number of other 
interesting properties; see {\it e.g.} \cite{B,Y} for more information. At about the same time,
building on Oseledec's theorem on Lyapunov exponents \cite{O}, 
Pesin \cite{P} and Ruelle \cite{R2} extended the uniform theory of hyperbolic
systems, also known as Axiom A theory, to an almost-everywhere theory in which
positive and negative Lyapunov exponents replace the uniform expansion and
contraction in Axiom A. The idea of an SRB measure was brought to this broader
setting and studied there by mostly Ledrappier and Young; see {\it e.g.}~\cite{LedY}.

\bigskip
\noindent {\bf The existence problem}

\smallskip
While the idea and relevant properties of SRB measures were shown
to make sense in this larger setting, existence was not guaranteed. 
Indeed for an attractor outside of the Axiom A category, no matter how chaotic it appears, there is, to this day, no general theory that will tell us whether or not it has
an SRB measure.

Here is where the difficulty lies: By definition, an Axiom A attractor has 
well-separated expanding and contracting directions that are
invariant under the dynamics, so that  tangent vectors in
 expanding directions are guaranteed some amount of growth with every iterate. 
In general, an attractor that appears chaotic to the eye 
must expand somewhere; this is how instabilities
are created. But since volume is decreased, there must also be
directions that are compressed. Without further assumptions, for most points
$x$ and tangent vectors $v$, 
$\|Df_x^nv\|$ will sometimes grow and sometimes shrink as a function of $n$.
To prove the existence of an SRB measure, one must show that
{\it on balance}, $\|Df_x^nv\|$ grows exponentially for certain coherent 
families of tangent vectors. The absence of cancellations between
expansion and contraction is what sets Axiom A attractors apart from 
general chaotic attractors.

%%%%%%%%%%%%%%%%%%%%%%%%%%%%%%%%%%%
\subsection{Some recent results on rank-one attractors} 

This subsection reviews some work by Wang and Young \cite{WY1, WY4, WY5}
on a class of strange attractors. These attractors have a single
direction of instability and strong contraction in all complementary
directions. Among systems without {\it a priori} separation of expanding
and contracting directions (or invariant cones), this is the only class
to date for which progress has been made on the existence of SRB
measures.

The idea is as follows: One embeds the systems of interest in a larger collection,
letting $b$ denote an upper bound on their contraction in all but one of
the directions (more precisely the second largest singular value of $Df_x$).
One then lets $b \to 0$ in what is called the {\it singular
limit}. {\it If } this operation results in
 a family of well defined 1D maps, and {\it if } some of 
these 1D maps carry strong enough expansion, then one can try to conclude that 
for small but positive $b$, some of the systems have SRB measures. 
Obviously, this scheme is relevant only for attractors that have a 1D
character to begin with. For these attractors, what is exploited here is the fact 
that 1D objects,
namely those in the singular limit, are more tractable than the original 
$n$-dimensional maps.

Since it is not illuminating to include all technical details in 
a review such as this one, we  refer the reader to \cite{WY4}, Sect.~1, for a formal statement, giving only enough information here to convey the flavor of the main result:

\bigskip
Let $M=I \times D_{n-1}$ where $I$ is either a finite interval or the circle 
$\mathbb S^1$ and $D_{n-1}$ is the closed unit disk in ${\mathbb R}^{n-1}$, 
$n \ge 2$. Points in $M$ are denoted by $(x,y)$ where $x \in I$ and 
$y=(y^1, \cdots, y^{n-1}) \in D_{n-1}$, and $I$ is sometimes identified with $I \times
\{(0,\cdots,0)\}$. Given $F: M \to I$, we associate 
two auxiliary maps:
\begin{align*}
F^\sharp : M \to M  & \quad  {\rm where} \quad 
F^\sharp = (F, 0, \cdots, 0)\ , \\
f: I \to I & \quad {\rm where} \quad f = F|_{I\times\{(0,\cdots,0)\}}\ . 
\end{align*}
We need to explain one more terminology: There is a well known class of
1D maps called {\it Misiurewicz maps} \cite{M}. Roughly speaking, a map $f$ is
in this class if it is $C^2$, piecewise monotonic with nondegenerate
critical points, and satisfies the following conditions: (i) it is
expanding away from $C=\{f'=0\}$, and (ii) the forward orbit of every $\hat
x \in C$ is trapped in an expanding invariant set (bounded away from $C$).
Maps in this class are known to have positive Lyapunov
exponents Lebesgue-a.e.

\begin{theorem}[\cite{WY4}]
Let $F_a : M \to I$ be a 1-parameter family of $C^3$ maps with the
following properties:

\smallskip
{\rm (C1)} there exists $a^*$ such that $f_{a^*}$ is a Misiurewicz map;

{\rm (C2)} $a \mapsto f_a$ satisfies a transversality condition at $a=a^*$ and 
$\hat x \in C(f_{a^*})$;

{\rm (C3)} for every $\hat x \in C(f_{a^*})$, there exists $j$ such that
$\partial_{y^j} F_{a^*}(\hat x, 0) \ne 0\ $.

\smallskip
\noindent Then there exists $b>0$ (depending on $\{F_a\}$) such that if
$T_a: M \to M$ is a family of $C^3$ embeddings of $M$ into itself with
$\|T_a - F_a^\sharp\|_{C^3} < b$, then there is a
positive measure set $\Delta$ in $a$-space such that for all $a \in \Delta$, 
$T_a$ admits an SRB measure.
\end{theorem}

That $\|T_a - F_a^\sharp\|_{C^3}$ must be sufficiently small is 
the {\it rank one} condition discussed above, and (C1) is where we require
the singular limit maps to have sufficient expansion. 
Notice that (C1)--(C3) all pertain to behavior at or near $a=a^*$.
The set $\Delta$ will also be in the vicinity of this parameter.
(C2) guarantees that one can bring about changes effectively by 
tuning the parameter
$a$, and (C3) is a nondegeneracy condition at the critical points.

\medskip
\noindent {\bf Remark.}  The existence of SRB measures is asserted for a
positive measure set of parameters, and not for, say, an entire interval
of $a$. This is a reflection of reality rather than a weakness of the
result: there are parameters $a$ arbitrarily near $\Delta$ for which
$T_a$ has sinks.  In a situation such as this one where chaotic and
non-chaotic regimes coexist in close proximity of one another, it is
impossible to say for certain if any given map has an SRB measure. One
can conclude, at best, that nearby maps have SRB measures ``with
positive probability".

\medskip
Theorem 1 was preceded by the corresponding result for the H\'enon family
\begin{equation}
\label{henon}
T_{a,b} : (x,y) \mapsto (1-ax^2+y, bx)\ .
\end{equation}
The existence of SRB measures for parameters near $a^*=2$ and $b\ll1$
was proved in \cite{BY} building on results from \cite{BC}. This is the first time 
the existence of SRB measures
was proved for genuinely nonuniformly hyperbolic attractors. 
Even though \cite{BC} is exclusively about Eq.~(\ref{henon}), the techniques 
developed there were instrumental in the proof of Theorem 1. 

\medskip

Returning to the setting of Theorem 1, let us call $a \in \Delta$ a
``good parameter" and $T=T_a$ a ``good map". The following two
properties of these maps are directly relevant to us. They were proved
under the following additional assumption on $M$:   
\begin{equation}\tag{$\star$}
  |\det(DT_a)| \sim b^{n-1}\ .
\end{equation}
\begin{itemize}
\item[(1)] Lebesgue-a.e. $z \in M$ is contained in $W^s(\xi)$ where $\xi$ is typical
with respect to an ergodic SRB measure (in general, there may be more than one such measure). It follows that $\Lambda_{\max}(z)>0$ 
for Lebesgue-a.e. $z \in M$. 

\item[(2)] Another condition on $f_{a^*}$ (Lyapunov exponent $> \log 2$ 
and $f_{a^*}^N$ mapping every interval of monotonicity to all of $I$ for some $N$)
implies the uniqueness of SRB measure. This in turn implies $\Lambda_{\max}$
is constant a.e. in $M$.
\end{itemize}

These and a number of other results for ``good maps" were proved in
\cite{WY5}.  We mention one that is not used here but sheds light on the
statistical properties of these attractors: For an SRB measure $\mu$ for
which $(T,\mu)$ is mixing, the system has exponential decay of
correlations for Lipschitz observables, {\it i.e.}, there exists $\tau
\in (0,1)$ such that for all Lipschitz $\varphi, \psi$, there exists
$C=C(\varphi, \psi)$ such that for all $n \ge 1$,
$$
\left| \int (\varphi \circ T^n) \psi~ d\mu - \int \varphi~ d\mu \int \psi~ d\mu \right|
\le C \tau^n\ .
$$

%%%%%%%%%%%%%%%%%%%%%%%%%%%%%%%%%%%%%
\subsection{Application to kicked oscillators}

We now return to the model introduced in Sect.~1.1 and explain how this system 
can be fitted into the framework of the last subsection. (See \cite{WY2} for details.)
We fix $\sigma, \lambda, A$, and allow $\tau$ to vary. Writing $\tau = k+a$
where $k=[\tau]$, the integer part of $\tau$, we let $T_{k,a}=\Psi_\tau$.
For each fixed $k \in \mathbb Z^+$, we view $\{T_a=T_{k,a}, a \in [0,1)\}$ as the family of
interest, and discuss if and when the conditions of Theorem 1 will hold
for this family.   

Here, the singular limit maps $F_a$ are well defined. In fact, they are
 the first components of $\lim_{k \to \infty}T_{k,a}$,
{\it i.e.},
$$
F_a(\theta, y) = \theta + a + \frac{\sigma}{\lambda} \cdot (y +
A \sin(2\pi \theta))\ ,
$$
and the restriction of $F_a$ to ${\mathbb S}^1$ is 
\begin{equation} \label{1d}
f_a: {\mathbb S}^1 \to {\mathbb S}^1, \qquad f_a(\theta) = \theta + a +
\frac{\sigma}{\lambda}A \sin(2\pi \theta), \qquad a \in [0,1)\ .
\end{equation}

Notice immediately that the range of applicability of Theorem 1 is limited to 
$\lambda \tau$ relatively large. This is because
$\|T_{k,a} - F_a^{\sharp}\|_{C^3} = {\cal O}(b)$ where 
$b=e^{-\lambda k}$ is required to be very small. 
For a given unforced system, 
where the amount of damping
$\lambda$ is fixed, this means the kicks must be applied sufficiently far apart
in time.

We comment on  (C1), which along with the rank one condition above
are the core assumptions for this theorem. For our purposes let us assume
$f_a$ satisfies the Misiurewicz condition if some iterate of $f_a$ sends 
its two critical points $c_1$ and $c_2$ 
into an unstable periodic orbit or an expanding invariant Cantor set. 
First, such Cantor sets are readily
available for medium size values of $\frac{\sigma}{\lambda}$ such as 
$\frac{\sigma}{\lambda} \ge 1$, and unstable periodic orbits start to exist for
somewhat smaller values of $\frac{\sigma}{\lambda}$. Suppose for some
parameter value $a$ that 
the forward orbit of $c_1$ is contained in an expanding invariant set $K$. 
As we vary $a$, both the orbit of $c_1$ and $K$ will move with $a$.
Condition (C2), assuming it holds, implies that for all large enough 
$k \in {\mathbb Z}^+$, 
$f^k_a(c_1)$ moves faster, {\it i.e.}, 
the path traced out by $a \mapsto f^k_a(c_1(a))$ cuts across $K$ as
though the latter was stationary. When $K$ is a Cantor set, this guarantees that  
 $f^k_a(c_1) \in K$ for an uncountable number of $a$'s. 
 By symmetry, when that happens to 
$f_a(c_1)$, the same is automatically true for $f_a(c_2)$. 
The larger $\frac{\sigma}{\lambda}$, the denser these Cantor sets are
in $\mathbb S^1$,
and the denser the set of parameters $a$ that can be taken to be $a^*$
in (C1).  

The checking of (C2), (C3), and ($\star$) are straightforward. 
Conditions for the uniqueness of SRB measures 
require that $\frac{\sigma}{\lambda}$ be a little larger. 

To summarize, the results in the last subsection imply that for 
$\frac{\sigma}{\lambda} \ge1$ (or even smaller),
for all large enough $k$, there exist positive measure sets $\Delta_k$ 
such that for $a \in \Delta_k$ and $\tau=k+a$,
$\Psi_\tau$ is a ``good map" in the sense of the last subsection.
In particular, $\Psi_\tau$ has an SRB measure. We conclude also that
$\Lambda_{\max}(z)$ is well defined and $>0$ for 
Lebesgue-a.e. $z \in U$. To ensure that  $\Lambda_{\max}(z) =$
constant a.e. (so there are no open
squares in Fig. 2) one needs to take $\frac{\sigma}{\lambda}$ a little larger.
It is in fact not hard to see that $\Delta_k \approx \Delta_{k'}$ 
for $k \ne k'$ when both are sufficiently large, so that in this parameter 
range, the set of $\tau$ for which 
the properties above are enjoyed by $\Psi_\tau$ 
is roughly periodic with period $1$.

\bigskip
\noindent {\bf Remarks on analytic results for chaotic systems}

\smallskip
Theorem 1 is a perturbative result. As is generally the case with 
perturbative proofs, the sizes of the perturbations (such as $b$) are 
hard to control. Consequently, applicability of Theorem 1 is limited to regimes with
very strong contraction. The results reported in Sect. 3.2, however, 
are the only rigorous results available at the present time. Techniques for analyzing
maps in  parameter ranges such as those in Fig.~2 are lacking and currently
quite far out of reach. 

The situation here is a reflection of the general state of affairs: 
Due to the cancellations discussed
at the end of Sect.~3.1, rigorous results for the large-time
behavior of chaotic dynamical systems tend to be challenging.

When results such as Theorem 1 are available, however, they -- and the ideas
behind them -- often shed light on situations that are technically
beyond their range of
applicability.  
Our example here is a good illustration of that:
 Fig.~2(d)-(f) show that as $\frac{\sigma}{\lambda}A$ increases, positive Lyapunov exponents become more abundant among the parameters
tested, interspersed with occasional sinks. This is in agreement 
with the dynamical picture suggested by Theorem 1 even though with
$\lambda \tau \in [0.5, 1.5]$, the contraction can hardly be considered strong.

%%%%%%%%%%%%%%%%%%%%%%%%%%%%%%%%%%%%%
%%%%%%%%%%%%%%%%%%%%%%%%%%%%%%%%%%%%%
\section{Generalizations}

In Sections 1, 2, and 3.3, we have focused on a concrete model.
We now generalize this example in two different ways:
\begin{itemize}
\item[$\bullet$] the unforced equation in Eq.~(\ref{model1}) is replaced
by an arbitrary limit cycle;
\item[$\bullet$] the specific kick in Eq.~(\ref{model1}) is replaced by an
arbitrary kick.
\end{itemize}
More precisely, we consider a smooth flow $\Phi_t$ on a finite
dimensional Riemannian manifold (which can be ${\mathbb R}^n$), and let
$\gamma$ be a {\it hyperbolic limit cycle}, {\it i.e.}, $\gamma$ is a
periodic orbit of period $p$, and for any $x \in \gamma$, all
eigenvalues for $D\Phi_p(x)$ are $<1$ aside from that in the flow
direction.  The {\it basin of attraction} of $\gamma$ is the set ${\cal
  B} := \{x \in M: \Phi_t(x) \to \gamma$ as $t \to \infty\}$. We
continue to consider forcing in the form of kicks, and assume for
simplicity that the kicks are defined by a  smooth embedding
$\kappa : M \to M$
(as would be the case if, for example, $\kappa$ represents the result of
a forcing defined by  $\dot z = G(z,t)$
where the vector field $G(z,t)$ is nonzero, or ``on,'' for only a short
time). As before, the kicks are applied periodically at times $0, \tau,
2\tau, 3\tau, \dots$, and the time evolution of the kicked system is
given by $\Psi_\tau = \Phi_\tau \circ \kappa$.

A new issue that arises in this generality is that one may not be able to
isolate the phenomenon, that is to say, the kicking may cause the limit cycle
to interact with dynamical structures nearby. Our discussion below is limited
to the case where this does {\it not} happen, {\it i.e.}, we assume 
there is an open set $U$ with $\gamma \subset U \subset {\cal B}$ 
such that $\kappa(U) \subset {\cal B}$ and $\Phi_\tau(\kappa(U)) \subset U$, 
and define $\Gamma = \cap_n \Psi_\tau^n({U})$ to be the {\it attractor} of 
the kicked system as before. 

We further limit the scope of our discussion in the following two ways:
(i) Kicks that are too weak will not be considered; such kicks produce
invariant curves and sinks for the same reasons given in Sect. 2.1, and
there is no need to discuss them further.  (ii) We consider only regimes
that exhibit a substantial contraction during the relaxation period,
brought about by long enough kick intervals that permit the ``shear" to act. 
As we will see, this is a more tractable situation. Rigorous results can
be formulated -- and we will indicate what is involved -- but will focus
primarily on ideas. Precise formulations of results in this generality
(see \cite{WY3}) are unfortunately not as illuminating as the phenomena
behind them.

\bigskip
\noindent {\bf The geometry of folding:  kicks and
the strong stable foliation}

\smallskip
As we will show, key to understanding the effects of kicks is the geometric
relation between the kick and the strong stable foliation associated with
 the limit cycle of
the unforced system. For $x \in \gamma$, we define the {\it strong
stable manifold} of $\Phi_t$ at $x$, denoted 
$W^{ss}(x) = W^{ss}_{\Phi_t}(x)$, to be the set 
$W^{ss}(x)=\{y \in M: d(\Phi_t(y), \Phi_t(x)) \to 0$
as $t \to \infty \}$; the distance between $\Phi_t(x)$ and
$\Phi_t(y)$ in fact decreases exponentially; see {\it e.g.} \cite{HPS}.
(This stable manifold is for the flow $\Phi_t$, not to be confused with that
for the kicked map $\Psi_\tau$ in Prop.~2.1.)       
Some basic          
properties of these manifolds are: (i) $W^{ss}(x)$ is a
codimension one submanifold transversal to $\gamma$ and meets
$\gamma$ at exactly one point, namely $x$; (ii)
$\Phi_t(W^{ss}(x)) = W^{ss}(\Phi_t(x))$, and in particular, if
the period of $\gamma$ is $p$, then
$\Phi_p(W^{ss}(x))=W^{ss}(x)$; and (iii) the collection
$\{W^{ss}(x), x \in \gamma\}$ foliates the basin of attraction of $\gamma$.
An example of a $W^{ss}$-foliation for a limit cycle is shown in Fig.~4. 

%%%%%%%%%%%%%%%%%%%%%%%%%%%%%%%%
%\ 
\begin{figure}
  \begin{center}
    \includegraphics*[scale=1.2,bb=80 0 373 105]{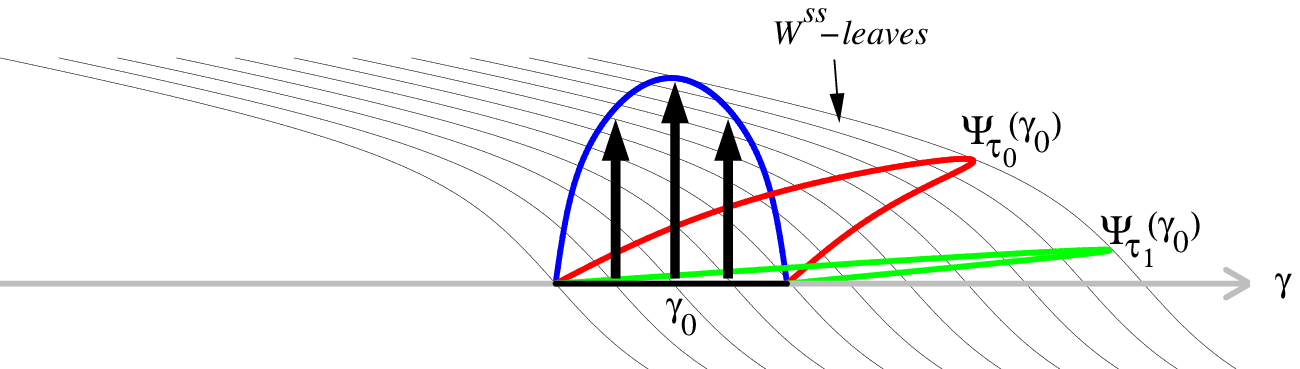}
  \end{center}

%\vskip 2in
  \caption{Geometry of folding in relation to the $W^{ss}$-foliation.
  Images of $\Psi_{\tau_1}(\gamma_0)$ and $\Psi_{\tau_2}(\gamma_0)$
  for $\tau_1< \tau_2$, both multiples of the period of the limit cycle $\gamma$, 
  are shown.}
\end{figure}  
%\bigskip  
%%%%%%%%%%%%%%%%%%%%%%%%%%%%%%%%

Fig.~4 shows the image of a segment $\gamma_0$ of $\gamma$ under
$\Psi_\tau$. For illustration purposes, we assume $\gamma_0$ is kicked
upward with its end points held fixed, and assume $\tau=np$ for some $n
\in {\mathbb Z}^+$ where $p$ is the period of the cycle.  Since
$\Phi_{np}$ leaves each $W^{ss}$-manifold invariant, we may imagine that
during relaxation, the flow ``slides'' each point of the curve
$\kappa(\gamma_0)$ back toward $\gamma$ along $W^{ss}$-leaves; the
larger $n$ is, {\it i.e.}, the more times it laps around,
the farther down the point slides. In the      
situation depicted, the folding is quite evident.  If $\tau$ is not an
integer multiple of $p$, then $\Phi_\tau$ carries each $W^{ss}$-manifold
to another $W^{ss}$-manifold. Writing $\tau = np + a$ where $a \in [0,
  p)$, we can think of the action of $\Phi_\tau$ as first sliding along the 
  $W^{ss}$-manifold by an amount corresponding to $\Phi_{np}$ and then
  flowing forward for time $a$.  
  
The picture in Fig.~4 gives insight into what types of kicks are likely to 
produce chaos. The following observations are intended to be informal
but intuitively clear:

\begin{itemize}

\item[(i)] Kicks directed along $W^{ss}$-leaves or in directions roughly
  parallel to $W^{ss}$-leaves are not effective in producing chaos, nor
  are kicks that essentially carry one $W^{ss}$-leaf to another in an 
  order-preserving fashion. For such kicks, $\Psi_\tau$ essentially permutes
$W^{ss}$-leaves, and $\kappa$ has to overcome the contraction within
  individual leaves to create chaotic behavior. (This cannot happen in 2D.)

\item[(ii)] The stretch-and-fold mechanism for producing chaos remains
  valid: the more $\Psi_\tau(\gamma)$ is folded, the more chaotic the
  system is likely to be, {\it i.e.}, the intuition is as in Fig.~1.
  What is different here is that unlike our earlier example, where the
  propensity for shear-induced chaos is determined entirely by
  parameters in the unforced equation, namely $\sigma$ and $\lambda$, we
  see in this more general setting that it matters how the kick is
  applied.  {\it It is the geometry of the action of $\kappa$ on the
    limit cycle $\gamma$ in relation to the strong stable foliation
    $W^{ss}$ that determines the stability or chaos of the kicked
    system.}

\item[(iii)] The case of stronger contraction is more tractable
  mathematically for the following reason: When the contraction in
  $\Psi_\tau$ is weak, as with $\tau=\tau_1$ in Fig.~4, one has to deal
  with the cumulative effects of multiple kicks, which are difficult to
  treat.  When the image $\Psi_{\tau}(\gamma_0)$ is pressed more
  strongly against $\gamma$, as in the case of $\tau=\tau_2 > \tau_1$,
  cumulative effects of consecutive kicks are lessened.

\end{itemize}
%% Some of the ideas above will be illustrated in the
%% examples to follow.

We illustrate some of the ideas above in the examples below.

\bigskip
\noindent {\bf Linear shear-flow examples}

\medskip
\noindent {\bf The 2D system in Eq.~(\ref{model1}):} The ideas in
this section 
can be seen as an abstraction of those discussed earlier. To understand that, we compute the 
$W^{ss}$-leaves of the unforced equation in Eq.~(\ref{model1}), and find them to be straight lines having slopes 
$-\frac{\lambda}{\sigma}$. We concluded earlier that
given $A$, the larger $\frac{\sigma}{\lambda}$, {\it i.e.}, the
smaller the angle between the $W^{ss}$-leaves and $\gamma$, the more chaotic 
the system is likely to be. Item (ii) above corroborates this conclusion: 
Given that we kick perpendicularly to the limit cycle (as is done in Eq.~(\ref{model1})),
 and points in $\gamma$ are kicked to a given height, 
the more ``horizontal" the $W^{ss}$-leaves, the farther the points in $\kappa(\gamma)$ 
will slide when we bring them back to $\gamma$. In other words, 
the $\frac{\sigma}{\lambda}$ part of the ratio
from earlier is encoded into the geometry of the $W^{ss}$-foliation ---
provided that we kick perpendicularly to the cycle. 

\medskip
\noindent {\bf Generalization to $n$ dimensions:} The $n$-dimensional analog of
Eq.~(\ref{model1}) with a more general forcing is
\begin{equation}
\begin{array}{ccl}
\dot{\theta} &=& 1 + \sigma \cdot {\bf y} \ ,\\[1ex]
\dot{\bf y} &=& -\Lambda {\bf y} + A H(\theta) {\bf v}(\theta) 
\cdot\sum_{n=-\infty}^\infty \delta(t-n\tau)\ ,\\
\end{array}
\label{model2}
\end{equation}
where $\theta \in {\mathbb S}^1$, ${\bf y} \in {\mathbb R}^{n-1}$,
$\sigma\in\R^{n-1}$ is nonzero,  and
$\Lambda$ is an $(n-1)\times (n-1)$ matrix all of whose eigenvalues have
strictly positive real parts. For simplicity, we assume the kicks are
perpendicular to the limit cycle $\{{\bf y}=0\}$, and to facilitate the
discussion, we have separated the following aspects of the kick
function: its amplitude is $A$, variation in $\theta$ is $H(\theta)$,
and the direction of the kick is ${\bf v}(\theta) \in {\mathbb
  S}^{n-2}$. As a further simplification, let us assume 
  ${\bf v}(\theta) \equiv {\bf v} \in 
  {\mathbb R}^{n-1}$, {\it i.e.}, it is a fixed vector.
  
 A computation shows that the $W^{ss}$-manifolds of the
unforced equation are given by  $$W^{ss}(\theta_0,{\bf 0}) = 
\{(\theta, {\bf y}) : \theta =
\theta_0 -\sigma^T{\bf \Lambda}^{-1}{\bf y}\},$$
 {\it i.e.}, they are
  hyperplanes orthogonal to the covector
  $(1,\sigma^T{\bf\Lambda}^{-1})$~.  
As noted in item (i) above, kick
  components orthogonal to $(1,\sigma^T{\bf \Lambda}^{-1})$ are
  ``dissipated'' and do not have much effect.
If $H \equiv $ constant, then $\Psi_{\tau}$
simply permutes the $W^{ss}$-planes and again no chaotic behavior will ensue.
To produce horseshoes and strange attractors, a sufficient amount of
{\it variation} in $\theta$ for $\Psi_\tau(\gamma)$ is needed as noted in
item (ii) above; that variation must come from $H$. An analysis similar to
that in Eq.~(\ref{time-t map}), Sect.~1.2, tells us that for large $\tau$, the
amount by which the kick is magnified in the $\theta$-direction is $
\approx A~H(\theta)~\sigma^T {\bf \Lambda}^{-1} {\bf v}$~.
We remark that the variation in $H$ is far more
 important than its mean value, which need not be $0$.  
   
 Finally, given $H \not \equiv $ constant, to maximize the variation of
 $\Psi_\tau(\gamma)$ in $\theta$ for large $\tau$, the discussion above
 suggests kicking  in a
direction ${\bf v}$ that maximizes $\sigma^T {\bf \Lambda}^{-1}{\bf v}$.
Under the conditions above, this direction is unique and is given by
${\bf v} = \pm({\bf \Lambda^T})^{-1}\sigma/|({\bf \Lambda^T})^{-1}\sigma|$.
Notice that this need not be the direction with the least damping or the direction
with maximal shear, but one that optimizes the combined effect of both.

\bigskip
%\noindent {\bf Remarks on analytic proofs}
\noindent {\bf On analytic proofs}

\smallskip
When $\Psi_\tau$ contracts strongly enough, the system falls into the
{\it rank one} category as defined in Sect.~3.2, independently of the
dimension of the phase space. This comes from the fact that $\Phi_t$
carries all points back to $\gamma$, which is a one-dimensional object.

If the folding (as described in item (ii) above) is significant enough
for the amount of contraction present, then horseshoes can be shown to
exist. Proving the existence of horseshoes with one unstable direction
is generally not very difficult, and not a great deal of contraction is
needed.

To prove the existence of strange attractors or SRB measures, 
the results in Sect.~3.2 are as applicable here as in our 2D linear example.
It is proved in \cite{WY3} that the periodic kicking of arbitrary
limit cycles fits the general framework of Theorem 1, in the sense
that as the time between kicks tends to infinity,
 singular limit maps are well defined. They are given by
$f_a:\gamma \to\gamma$,
$a \in [0,p)$, where
$$
f_a(x) \ := \ \lim_{n \to
  \infty}\Phi_{np+a}(\kappa(x))\qquad\mbox{ for all }
x\in\gamma\ .
$$
 From our earlier discussion of sliding along $W^{ss}$-leaves, 
it is not hard to see that $f_a(x)$ is, in fact, the unique point $y \in \gamma$ such
that $\kappa(x) \in W^{ss}(y)$. Whether  (C1)--(C3) hold 
depends on the system in question and hinges mostly on (C1), which 
usually holds when the variation is large enough. As always, these conditions 
need to be verified from example to example.

%%%%%%%%%%%%%%%%%%%%%%%%%%%%%%%%
\section*{Related Results and Outlook}

We have reviewed a set of results on the periodic kicking of limit
cycles. The main message is that the effect of
the kick  can be magnified by the 
underlying shear in the unforced system to create an unexpected amount of 
dynamical complexity. It is an example of a phenomenon known as
 {\it shear-induced chaos}. 
 
A similar geometric mechanism is used to prove the existence of strange
attractors in (a) certain examples of slow-fast systems \cite{GWY}; (b)
periodic kicking of systems undergoing supercritical Hopf bifurcations
(see \cite{WY3} for details, and \cite{LWY} for results applicable to
evolutionary PDEs); 
and (c) periodic forcing of near-homoclinic loops \cite{WO}. 
See also \cite{OS}. All of these results pertain to strong-contraction regimes;
proofs are perturbative and rely on the theory of rank one attractors
reviewed in Sect.~3.2.

A welcome extension of the results reviewed here is to remove the
strong-contraction assumption for strange attractors, but this is likely
to be challenging: one has to either develop non-perturbative techniques
or go about the problem in a less direct way.

{\it Random forcing} is a future direction we believe to be both interesting
and promising. Numerical studies of
Poisson and white-noise forcing have been carried out \cite{LY}.
Phenomena similar to those in Sect.~1 are observed when the kick term in
Eq.~(\ref{model1}) is replaced by a term of the form
$A\sin(2\pi\theta)~dB_t$ where $B_t$ is standard Brownian motion, {\it
  i.e.}, the setup is a stochastic differential equation.  With
stochastic forcing, phase space geometry is messier, but $\lmax$ depends
continuously on parameters. The absence of wild fluctuations between
positive and negative values of $\lmax$ (corresponding respectively to 
strange attractors and sinks in the periodic case) gives hope
to the idea that the analysis for
stochastic forcing may be more tractable than that for periodically forced systems.

%% Positive Lyapunov exponents are seen even under weak forcing in
%% mildly contracted regimes, due possibly to large deviations.

%%%%%%%%%%%%%%%%%%%%%%%%%%%%%%%%

\end{document}